\documentclass[12pt]{amsart}
\usepackage{amscd,graphicx,epsfig}

\theoremstyle{definition}

\def \ph{\varphi}

\def \diag{\operatorname {diag}}

\def \refeq#1{equation (\ref{#1})}
\def \ra{\rightarrow}

\def \hom{\mbox{\rm Hom}}
\def \ie{\hbox{\it i.e.}}

\def \mcom{,\cdots,}
\def \k{\mbox{$\mathbb K$}}

\def \C{\mbox{$\mathbb C$}}
\def \Z{\mbox{$\mathbb Z$}}

\def\zt{\mbox{$\Z_2$}}

\def\ad{\operatorname{ad}}

\def\inv{^{-1}}
\def\d{d}

\def\im{\operatorname{im}}
\def\A{\mbox{$\mathcal A$}}

\def\W{W}

\def\coder{\operatorname{Coder}}

\def\linf{\mbox{$L_\infty$}}
\def\and{\mbox{ \rm and }}

\def\psa#1#2{\psi^{#1}_{#2}}

\def\phd#1#2{\ph^{#1}_{#2}}

\def\inv{^{-1}}

\def\ho{\text{\it ho}}

\def\aut{\operatorname{Aut}}
\def\dinfty{\mbox{$d^\infty$}}
\def\dinf{\mbox{$d^\text{inf}$}}
\def\P{\mathbb P}
\author{Marilyn Daily}
\address{Max Planck Institut f\"ur Gravitationsphysik\\
14476 Golm, Germany} \email{Marilyn.Daily@aei.mpg.de}
\author{Alice Fialowski}
\address{E\"otv\"os Lor\'and University\\
P\'azm\'any P\'eter s\'et\'any 1/C\\
H-1117 Budapest, Hungary} \email{fialowsk@cs.elte.hu}
\author{Michael Penkava}
\address{University of Wisconsin\\
Eau Claire, WI 54702-4004}
\email{penkavmr@uwec.edu}
\thanks{Research of these authors was partially supported by grants
from the National Science Foundation, Hungarian grants OTKA T043641 and
T043034, grants from the University of Wisconsin-Eau Claire, and
Czech 
grant ME 603.}
\title[3-dimensional \linf\ algebras]
{Comparison of 3-dimensional \Z-graded and \zt-graded \linf\ algebras}
\begin{document}
\setlength{\multlinegap}{0pt}
\begin{abstract}
This article explores \linf\ structures
-- also known as \lq strongly homotopy Lie algebras\rq\ --
on 3-dimensional vector spaces with both \Z- and \zt-gradings.
Since the \Z-graded \linf\ algebras are special cases of \zt-graded
algebras in the induced \zt-grading, there are
generally fewer \Z-graded \linf\ structures on a given space.
On the other hand, only degree zero automorphisms, rather than
just even automorphisms, are used to determine equivalence in a \Z-graded space.
We therefore find nontrivial examples in which the map from the \Z-graded moduli space
to the \zt-graded moduli space is bijective, injective but not surjective, or
surjective but not injective.
Additionally, we study how the codifferentials in the moduli
spaces deform into other nonequivalent codifferentials, which
gives each moduli space a sort of topology.
\end{abstract}
\date\today
\maketitle

\section{Introduction}Infinity algebras commonly arise in both \Z-gradings
and \zt-gradings.  For example, in the study of supersymmetry, a
\zt-grading is natural.  In topological or algebraic contexts,
though, a \Z-grading is often preferable.  In this article, we
explore the relationship between \linf\ structures with these two
different gradings.  Because both the \Z-graded and \zt-graded
\linf\ algebras on three-dimensional spaces have been thoroughly
studied in \cite{dal1,bfp1,fp3,fp6,fp7}, we will use \linf\
structures on these spaces to illustrate this relationship.

Every \Z-graded space has a natural \zt-grading, and since the
signs which arise in the computations depend only on this
\zt-grading, every \Z-graded infinity algebra induces a \zt-graded
infinity algebra.  Because of additional restrictions on the
structures in the \Z-graded case, not every \zt-graded infinity
algebra on a space will be a \Z-graded infinity algebra.  From a
naive point of view, it would seem that the \Z-graded infinity
algebras are just a subset of the \zt-graded infinity algebras.

However, the moduli space of infinity algebras consists of
equivalence classes of isomorphic structures, and the notion of
isomorphism depends on whether the structure is considered as
\Z-graded or \zt-graded. Since every isomorphism in the \Z-graded
category is also an isomorphism in the \zt-graded category, there
is a natural map from the moduli space of \Z-graded infinity
algebras to the moduli space of \zt-graded infinity algebras. The
purpose of this paper is to explore this map. We shall see that in
some cases the map is bijective, but in others the map may fail to
be injective or surjective.

In the \zt-graded category, there are only a few types of
three-dimensional vector spaces.  The dimension of a \zt-graded
space is denoted in the form $m|n$, where $m$ is the number of
even basis elements and $n$ is the number of odd basis elements.
Thus, three-dimensional \zt-graded algebras are all of dimension
$3|0$, $2|1$, $1|2$ or $0|3$. There are an infinite number of
possible three-dimensional \Z-graded spaces, which may be
identified with a triple $(k,m,n)$ of nondecreasing integers, representing
the \Z-graded degrees of the three basis vectors. Each \Z-graded space
also has a \zt-graded dimension, depending on which basis vectors are even and odd.
Most \Z-graded spaces have only a finite number of \linf-algebra
structures.  However, a few special cases where the integers are
consecutive (studied in \cite{dal1}) or where two of the integers
coincide (see \cite{dal2}) have a richer moduli space. We will
concentrate our study on these more interesting cases.

It can easily be 
shown that any nontrivial \linf\ structure on a vector space of type $0|3$ must be a
classical Lie algebra, and that there are no nontrivial
\linf\ structures on spaces of type $n|0$ (for arbitrary $n$), so we do not have to
consider spaces  of type $3|0$. Thus, the interesting examples are of
\Z-graded spaces of \zt-graded dimension $1|2$ or $2|1$, which allows us to draw
upon earlier work concerning \zt-graded spaces.
In \cite{fp6}, \linf\ structures on \zt-graded spaces of dimension
$1|2$ were completely classified, and 
\cite{bfp1} classified \linf\ structures on \zt-graded spaces of type $2|1$.

Since the \zt-grading determines the signs in the computation of
brackets, 
it is not surprising that \Z-graded spaces which are of
dimension $2|1$ in the induced \zt-grading have a completely
different bracket structure than those whose induced \zt-grading
is of dimension $1|2$. Accordingly, we will treat these two cases
separately.
It is interesting to note that
since the signs in the bracket computations are 
entirely determined by the induced \zt-grading,
the brackets of elements in a \Z-graded space 
can be calculated using exactly the same
formulas that are used in the \zt-graded case.

The structure of this paper is as follows. Sections 2 through 4
cover preliminary notions about the moduli space of \linf\
algebras and the deformation theory of these algebras. In section
\ref{sec1bar2},  we consider the most interesting \Z-graded
example of dimension $1|2$, the vector space of type $(-1,0,1)$,
in considerable detail.  In this particular example, we show that
the \Z-graded moduli space is the same as the \zt-graded moduli
space, down to the smallest detail that can be obtained from the
deformation picture.

In section \ref{sec2bar1}, we will consider \Z-graded spaces of
dimension $2|1$, where there are two interesting examples.
The moduli space of the \Z-graded vector space of type $(0,1,2)$
is relatively simple, and it maps injectively (but not
surjectively) to the \zt-graded moduli space.
For a \Z-graded vector space of type $(-2,-1,0)$, the map of the
moduli spaces of codifferentials of order less than three {\em is}
surjective. However, this map fails to be injective, due to a
symmetry arising from an action of the symmetric group $\Sigma_2$
on the \zt-graded space that does not appear in the \Z-graded
case.  Nevertheless, the structures of the moduli spaces are
surprisingly similar.

For each space that we study, we will first establish a
classification of the codifferentials of a fixed degree. Next, we
will use cohomology in order to determine the nonequivalent
extensions. These are the elements of our moduli space. Finally,
in order to understand how the space is glued together, we will
study versal deformations of each of our elements.
Since our primary aims are to assemble a picture of the moduli
space and to compare the \Z-graded and the \zt-graded moduli
spaces, we will provide computational details only when they are
essential. The calculations are given in more detail in
\cite{bfp1,fp6}.

We would like to thank Jim Stasheff for raising this question.

\section{The moduli space of \linf\ Algebras}\label{secmodspace}
Let $V$ be a graded vector space. If $V$ is \Z-graded, then the
\emph{desuspension} $W$ of $V$ is defined by $W_i=V_{i+1}$. For a
\zt-grading, this simply corresponds to reversing the parity of
each element, and the resulting \zt-graded
space is called the \emph{parity reversion}.
Originally, \linf\ algebras \cite{ls,lm} were
defined as structures on the exterior algebra $\bigwedge V$,
motivated by the definition of a Lie algebra, which is given by a
map $l_2:\bigwedge^2 V\ra V$. One can translate the definition of
an infinity algebra to the desuspended $W$ picture, where it becomes a
structure on the symmetric coalgebra $S(W)$. Although less natural,
it is simpler to express the notion of an \linf\ algebra in the
$W$ picture, because then an \linf\ algebra is simply a
codifferential on the (graded) symmetric coalgebra of $W$
\cite{sta4,ps2}. This is
the point of view we will take in this article.

A codifferential is a coderivation whose square is zero. For
\zt-graded \linf\ algebras, the codifferential is required to be
odd, while for \Z-graded \linf\ algebras, the codifferential is
required to have degree 1. (Note that it is also possible to work
with the suspension, rather than the desuspension of $V$, in which
case the codifferentials will have degree $-1$.)
Every \Z-graded vector space is
naturally \zt-graded, and coderivations with respect to the
\Z-grading are also coderivations with respect to the \zt-grading.
Since coderivations of degree 1 are odd, one can view \Z-graded
\linf\ algebras as special cases of \zt-graded \linf\ algebras.

If $d$ and $d'$ are codifferentials on $W$ and $W'$ resp., then a
morphism of \linf\ algebras is a coalgebra morphism $g:S(W')\ra
S(W)$ satisfying $gd'=dg$. When $W=W'$ and $g$ is an automorphism
of $S(W)$, then $d'=g\inv d g$, and we say that $d'\sim d$. Thus
$\aut(S(W))$ acts on the \linf\ structures, and we denote
$g^*(d)=g\inv d g$. The equivalence classes of this action are
called the \emph{moduli space} of \linf\ algebras on $W$. The
classification problem for \linf\ algebras is the study of the
structure of this moduli space.

Note that for \Z-graded \linf\ algebras, we require a morphism of
symmetric algebras to be a degree 0 map, while for \zt-graded
algebras, the map is required to be even.  Thus every morphism of
\Z-graded symmetric algebras is also a morphism of \zt-graded
symmetric algebras. From this it follows that two \linf\ algebras
which are equivalent as \Z-graded \linf\ structures are equivalent
as \zt-graded \linf\ algebras.  Thus there is a map from the
moduli space of \Z-graded \linf\ algebra structures to the moduli
space of \zt-graded \linf\ algebras. There are fewer \Z-graded
\linf\ structures on $W$, but at the same time, the notion of
equivalence is more restrictive, so it is not obvious whether the
map between these moduli spaces is injective or surjective.  As we
shall see, in general it is neither.

The main goal of this paper is to study the relationship between
the moduli space of \Z-graded \linf\ algebras on a
three-dimensional vector space $W=W_m\oplus W_{m+1}\oplus W_{m+2}$
and the corresponding moduli space of \zt-graded \linf\ algebras.
The \Z-graded \linf\ structures on such a vector space were
studied in \cite{dal1}, while \zt-graded \linf\ structures on a
three-dimensional space were studied in \cite{bfp1,fp3,fp6}.

Define $C^r=\hom(S^r(W),W)$ to be the space of $r$-cochains, that
is, the cochains of exterior degree $r-1$, and $C^r_s$ to be the
$r$-cochains of internal degree $s$. We will omit the adjectives
internal and exterior when the meaning is clear from the context.
A \Z-graded codifferential $d$ lies in
$C^\bullet_1=\prod_{r=1}^\infty C^r_1$, so we can express
$d=d_1+d_2+\cdots$, where $d_r\in C^r_1$. There is a natural
isomorphism between $C=\prod_{r=1}^\infty C^r$ and $\coder(S(W))$,
the space of graded coderivations of the symmetric algebra on $W$,
so that $C$ inherits the structure of a graded Lie algebra.
Moreover, $[C^k_m,C^l_n]\subseteq C^{k+l-1}_{m+n}$.

The condition that $d$ is a codifferential is that $[d,d]=0$. This
is equivalent to the relations
\begin{equation}\label{linf_rels}
\sum_{r+s=n+1}[d_r,d_s]=0 \ \ \forall\ n\in\mathbb N.
\end{equation}
If $k$ is the least integer such that $d_k\ne0$, then $d_k$ is
called the \emph{leading term} of $d$ and $k$ is called the
\emph{order} of $d$, \ie, the order of $d$ is the degree of the least
nonvanishing term in $d$.
Note that although $d_1$ was nonzero in many of the originally 
studied \linf\ structures, it is entirely possible for $d_1$ to
be zero ({\it e.g.}, when a classical graded Lie algebra is
viewed as an \linf\ structure).
Since $d_k$ is itself a codifferential, we can classify
all possible \linf\ structures by first identifying all
nonequivalent codifferentials $d_k$, and then finding
all of their nonequivalent extensions 
(to include those codifferentials which have higher order terms).

If we define $D:C\ra C$ by $D(\phi)=[d_k,\phi]$, then $D$ is a
codifferential on $C$; \ie, $D^2=0$.
This allows us to define the cohomology $H(d_k)$ as
$H(d_k)=\ker(D)/\im(D)$. Moreover, one can define the $n$-th cohomology
group  $$H^n(d_k)=\ker(D:C^n\ra C^{n+k-1})/\im(D:C^{n-k+1}\ra C^n).$$
We can re-express the relations
in \refeq{linf_rels} in the form
\begin{equation}\label{coho}
D(d_{n+1})=-\tfrac12\sum_{r=k+1}^n[d_r,d_{k+n+1-r}] \ \ \forall\
n\geq k.
\end{equation}
This formula is familiar to deformation theorists, and essentially
implies that \linf\ algebras are deformations of the
codifferential of the leading term.  In fact, one can show that if
the relations hold up to order $n$, then the term on the right
hand side in \refeq{coho} is always a cocycle with respect to the
coboundary operator $D$. A coderivation $d$ which satisfies the
relations \refeq{coho} for all $n\leq m$ is called an $L_m$
algebra.

We can interpret the relations as saying that an $L_m$ algebra can
be extended to an $L_{m+1}$ algebra structure precisely when the
cocycle on the right hand side of \refeq{coho} is a coboundary. We
will use this observation extensively in our construction of
equivalence classes of \linf\ structures.

An automorphism $\lambda$ of $S(W)$ is called \emph{linear} if
$$\lambda(w_1\cdots w_n)=\lambda(w_1)\cdots\lambda(w_n).$$
Note that this implies that $\lambda:W\ra W$ is an isomorphism.
Not every automorphism is linear; in fact, if an even (degree 0
for the \Z-graded case) element $\phi\in C^n$ for $n>1$ is
considered as a coderivation of $S(W)$, then
$\exp(\phi)=\sum_{i=0}^\infty\tfrac 1{i!}\phi^i$ is an
automorphism of $S(W)$.  Moreover, every automorphism $g$ can be
uniquely expressed in the form
$$g=\lambda\prod_{n=2}^\infty \exp(\phi_n),$$
where $\lambda$ is a linear automorphism and $\phi_n\in C^n$,
and every such expression yields a well defined automorphism. The
term $\lambda$ is called the \emph{linear} part of $g$, and the
term $\prod_{n=2}^\infty \exp(\phi_n)$ is called the \emph{higher
order} part of $g$. We can express the action of $g$ on
codifferentials by the formula
\begin{equation}\label{auto}
g^*=\left(\prod_{n=\infty}^2\exp(-\ad_{\phi_n})\right)\lambda^*,
\end{equation}
which means that the action of the higher order part of $g$ on $C$
can be computed in terms of brackets of coderivations.

If $d$ is a coderivation, then its automorphism group $\aut(d)$ is
the subgroup of $\aut(S(W))$ which leaves $d$ invariant. For a
codifferential $d_k$ of fixed degree, it is easy to see that its
automorphism group is generated by the linear automorphisms
$\lambda$ fixing $d_k$, together with the exponentials of the even
(degree 0) cocycles for the coboundary operator $D(d_k)$
determined by $d_k$.

If $d$ and $d'$ are two equivalent codifferentials, then their
leading terms are of the same degree $k$ and linearly equivalent;
\ie, $d'_k=\lambda^*(d_k)$ for some linear automorphism $\lambda$.
Thus the first step in classifying the equivalence classes of
codifferentials on $S(W)$ is to find the linear equivalence classes
of degree $k$ codifferentials, which we call the moduli space of
degree $k$ codifferentials.

If $d=d_k+d_l+\ho$, where $k<l$ and  $\ho$ refers to terms whose
degrees are higher than $l$, then $d_l$ is called the \emph{second
term} of $d$. It is immediate that $D(d_l)=0$, so $d_l$ is a
cocycle. A codifferential is said to be in  \emph{standard form}
if its second term, if any, is a nontrivial cocycle. If the second
term of a codifferential is a coboundary, say $d_l=-\tfrac 12
D(\gamma)$ for some $\gamma\in C^{l-k+1}$, then
$d'=\exp(-\ad_\gamma)d$  is an equivalent codifferential whose
second term has higher degree than $l$. Thus every codifferential
is equivalent to one in standard form, and the study of
nonequivalent extensions of $d_k$ can be reduced to the study of
nonequivalent extensions of codifferentials in standard form.

 If two extensions of $d_k$ are in standard form, then the
second term of each has the same degree $l$, and the cohomology
classes of these second terms are linearly equivalent. We say that
two cohomology classes $\bar\alpha$ and $\bar\beta$ are linearly
equivalent if there is a linear automorphism $\lambda$ of $d_k$
such that $\lambda^*(\bar\alpha)=\bar\beta$. Note that this action
is well defined because $\lambda^*$ leaves $d_k$ invariant.

If $d=d_k+d_l$ is a codifferential and $H^n(d_k)=0$ for all $n>l$,
then any extension of $d$ is equivalent to $d$. This observation
will be used extensively in this paper to simplify the
classification of \Z-graded codifferentials.

\section{Deformations of \linf\ Algebras}
For details of the theory of deformations of \linf\ algebras, see
\cite{fp1}. Here we recall the basic notions. Suppose that $d_k$ is a
codifferential, and $\{\delta_i\}$ is a complete set of
representatives for the cohomology with regard to $d_k$,
in other words, a basis of a preimage
of the cohomology (where $i$ is an index for the basis, not the degree of $\delta_i$).
Then the {\em universal infinitesimal deformation} $\dinf$ of
$d_k$ is given by
\begin{equation*}
\dinf=d_k+\delta_i t^i,
\end{equation*}
where the $t^i$'s are parameters. (Note that we use the Einstein
summation convention.)
The \emph{miniversal deformation} is the most general extension of
the universal infinitesimal deformation to a formal deformation,
and its primary importance is that it induces all possible formal
deformations. We will use the miniversal deformation to determine
the neighborhood of a given codifferential in the moduli space. In
\cite{fp1}, miniversal deformations were constructed for \linf\ algebras of \emph{finite
type}, which includes the finite dimensional case we are considering here.

To construct the miniversal deformation, we compute
$[\dinf,\dinf]$, and look at the higher order terms which appear.
Since the bracket must be a cocycle, it can be expressed as a sum
of coboundary terms and nontrivial cocycles.  The coboundaries may
be eliminated by adding higher order terms, so that eventually one
obtains a coderivation $\dinfty$ which has the property that its
self-bracket contains no coboundaries. However, the nontrivial
cocycles appearing in the bracket cannot be eliminated by this
process. Accordingly, we obtain some relations on the parameters
which must be satisfied if the bracket is to vanish.  These
relations are called the relations on the base of the versal
deformation, which is the algebra $\k[[t^1,t^2,\dots]]$ modulo the
ideal generated by these relations. (Here, $\k$ is the base field, which
we will assume later is just $\C$.)

In the process which we are describing, one can proceed level by
level, constructing at each stage an $n$-th order deformation, for
which the bracket vanishes up to order $n+1$, except for terms
involving the nontrivial cocycles, which give rise to  $n$-th
order relations. In the limit we obtain the miniversal deformation
$\dinfty$ and relations which are given as a formal power series.

The process described above potentially requires an infinite
number of steps, although, in practice, it often terminates after
a finite number of steps. It is often more practical to construct
the miniversal deformation  recursively as follows. Set
\begin{equation*}
\dinfty=d+\delta_i t^i +\gamma_i x^i,
\end{equation*}
where $\{\delta_i\}$ are representatives of the odd cohomology
classes, $\{\gamma_i\}$ are preimages of the even coboundaries,
and $x^i$ are formal power series in the parameters $t^i$. To
determine the coefficients $x^i$, we express
\begin{equation*}
[\dinfty,\dinfty]=\alpha_i r^i +\beta_i s^i +\tau_i y^i,
\end{equation*}
where $\{\alpha_i\}$ are representatives of the even cohomology
classes, $\{\beta_i\}$ is a basis of the even coboundaries, and
$\alpha_i,\beta_i,\tau_i$ are a basis of the even cochains. It
follows from \cite{fp1} that $s^i=0$ for all $i$, $r^i$ are formal
power series in the parameters $t^i$, and
$y^i=0\mod{(r^1,\dots)}$.  Moreover, the equations $s^i=0$ can be
solved to obtain the expressions for $x^i$ as formal power series
in the parameters $t^i$.

Suppose that we have constructed a miniversal deformation
\begin{equation}
\dinfty=d+\delta_it^i +\cdots.
\end{equation}
Suppose some set of values of the parameters $t^i$ satisfy the
relations on the base. The codifferential obtained by substituting
these values into the formula for the miniversal deformation is a
well-defined codifferential in the moduli space. We can think of
the values of the parameters which satisfy the relations as
determining a variety in the $(t^i)$ space, and the miniversal
deformation as a description of the elements in the moduli space
which are near to it.  Thus the miniversal deformation provides
information about how the moduli space is glued together.

A one parameter family $d(t)$ of deformations of $d_k$ is a curve
in the $(t^i)$ space such that $d(0)=d_k$.  It can happen that
there is a codifferential $d'$ such that $d(t)\sim d'$ whenever
$t\ne 0$. In this case, the one-parameter family of deformations
is called a {\em jump deformation}. Otherwise, $d(t)$ is called a
family of smooth deformations. It is interesting to note that jump
deformations are one way transformations; if $d$ has a jump
deformation to $d'$, then there is never any jump deformation in
the opposite direction.

Suppose that there is a smooth family $d(t)$ of deformations of $d$.
Then this family is said to factor through a jump deformation from
$d$ to $d'$ if there is a family $d'(t)$ of deformations of $d'$ such
that $d'(t)\sim d(t)$ whenever $t\ne0$. In this case, we do not consider
the deformations $d(t)$ to belong to the same \emph{stratum} of the moduli space
of $d$, which is determined by the smooth families of deformations of $d$
which do not factor through jump deformations.  We will not explain here
why the notion of stratum is well defined, but use this notion only for describing the moduli
space.   The pattern we will observe in this paper is that the moduli
space has a decomposition in terms of these strata, and that each stratum is an orbifold.
The smooth deformations which do not factor through jump deformations determine how a
stratum is glued together, and the jump deformations determine how the strata are glued
to each other.

\section{Notation}

It is customary when working with \zt-graded spaces to list the
even elements of a basis first, and we will follow this convention
even though it means that our basis elements will not be listed in
the consecutive order that is natural from the \Z-graded point of
view.
In order to express our codifferentials in a compact form, we will
use the following notation. Let $\{e_1,e_2,e_3\}$ be a fixed basis
of $W$. Then $\ph^{i_1 i_2 i_3}_k\in\hom(S^{i_1+i_2+i_3}(W),W)$ is
defined by
\begin{equation*}
\ph^{i_1 i_2 i_3}_k(e_1^{j_1}e_2^{j_2}e_2^{j_2}) = i_1!\ i_2!\
i_3!\ \delta^{i_1}_{j_1}\,\delta^{i_2}_{j_2}\,\delta^{i_3}_{j_3}\
e_k.
\end{equation*}
The factorials in the above definition give rise to simpler formulas
for the brackets of coderivations.
In order to make it easier to distinguish between  even and odd
coderivations, if $\ph^{ijk}_l$ is odd, we will denote it as
$\psi^{ijk}_l$ instead.
In Section \ref{secmodspace}, we defined $C^r=\hom(S^r(W),W)$ to
be the space of $r$-cochains, and $C^r_s$ to be the cochains of
internal degree $s$. Then $\ph^{i_1 i_2 i_3}_k\in
C^{i_1+i_2+i_3}_{s}$, where $s = |e_k| - (i_1 |e_1| + i_2 |e_2| +
i_3 |e_3|)$ is the difference of the degrees of the output and
input of $\ph^{i_1 i_2 i_3}_k$
(\ie, the degree of  $\ph^{i_1 i_2 i_3}_k$ as a map).

From now on, let us assume that the base field is $\C$.  Note that in our definition
of exponentials of coderivations, we implicitly have been assuming that the field
has characteristic zero.

\section{Codifferentials on a $1|2$-dimensional space} \label{sec1bar2}

If $W$ is of the form $W=W_{2x-1}\oplus W_{2x}\oplus W_{2x+1}$,
then as a \Z-graded space, it has type $(2x-1,2x,2x+1)$, but its
dimension is always $1|2$ as a \zt-graded space. In order to list
the even basis elements first, we set $|e_1|=2x$, $|e_2|=2x-1$ and
$|e_3|=2x+1$. Then we decompose the space of cochains as follows.
\begin{align*}
C^r_{2s}=&
\begin{cases}
\langle\phd{r,0,0}1,\phd{r-2,1,1}1,\phd{r-1,1,0}2,\phd{r-1,0,1}3\rangle,&
s=-x(r-1)\\
\langle\phd{r-1,0,1}2\rangle&s=-x(r-1)-1\\
\langle\phd{r-1,1,0}3\rangle&s=-x(r-1)+1
\end{cases}\\
C^{r}_{2s+1}=&
\begin{cases}
\langle\psa{r-1,1,0}1,\psa{r,0,0}3,\psa{r-2,1,1}3\rangle&s=-x(r-1)\\
\langle\psa{r-1,0,1}1,\psa{r,0,0}2,\psa{r-2,1,1}2\rangle&s=-x(r-1)-1
\end{cases}
\end{align*}
Note that $C^1_1$ has dimension 2 for all $x$. When $x=-1$,
$C^2_1$ has dimension 3. Otherwise, $C^r_1=0$ for all $r>1$,
unless $x=0$, in which case $C^r_1$ has dimension 3 for all $r>1$.
We study the  most important case, when $x=0$, but the case $x=1$
has some interesting features as well.

\subsection{Codifferentials on
$\mathbf{W=W_{-1}\oplus W_{0}\oplus W_{1}}$}\label{GoodExample}

The basis elements of $W=W_{-1}\oplus W_{0}\oplus W_{1}$ have
degrees $|e_1|=0$, $|e_2|=-1$ and $|e_3|=1$, and the space of
cochains has the following decomposition.
\begin{align*}
C^r_{-2}=&\langle\phd{r-1,0,1}2\rangle\\
C^r_{-1}=&\langle\psa{r-1,0,1}1,\psa{r,0,0}2,\psa{r-2,1,1}2\rangle\\
C^r_{0}=&\langle\phd{r,0,0}1,
\phd{r-2,1,1}1,\phd{r-1,1,0}2,\phd{r-1,0,1}3\rangle\\
C^r_1=&\langle\psa{r-1,1,0}1,\psa{r,0,0}3,\psa{r-2,1,1}3\rangle\\
C^r_2=&\langle\phd{r-1,1,0}3\rangle
\end{align*}
We shall refer to the bases above as the standard bases of these
spaces. Note that for $r=1$,  there are no terms with $r-2$ in
their superscripts.
Since a codifferential for a \Z-graded \linf\ algebra must have
degree $1$, any codifferential of order $k$ is of the form
\begin{equation*}
d=\psa{k-1,1,0}1a+\psa{k,0,0}3b+\psa{k-2,1,1}3c.
\end{equation*}
Since $[d,d]=\phd{2k-2,1,0}32b(ka-c),$ it follows that $d$ is a
codifferential precisely when $b=0$ or $c=ka$. 
We will now show that any such codifferential is equivalent to one
of the three types below:
\begin{equation}\label{codiffs1}
\begin{array}{lcl}
d_k(\lambda:\mu)&=&\psa{k-1,1,0}1\lambda+\psa{k-2,1,1}3\mu,\\
d_k^*&=&\psa{k,0,0}3,\\
d_k^\sharp&=&\psa{k-1,1,0}1+\psa{k,0,0}3+\psa{k-2,1,1}3k.
\end{array}
\end{equation}
If $b=0$, then we obtain the codifferential
\begin{math}
d_k(\lambda:\mu)=\psa{k-1,1,0}1\lambda+\psa{k-2,1,1}3\mu.
\end{math}
A diagonal linear automorphism  takes $d_k(\lambda:\mu)$ to a
multiple of itself, so that the codifferentials of this type are
parameterized by $(\lambda:\mu)\in\P^1(\C)$. If $b\ne0$, we must
have
\begin{equation*}
d=\psa{k-1,1,0}1\lambda+\psa{k,0,0}3\mu+\psa{k-2,1,1}3k\lambda.
\end{equation*}
This time, if $g=\diag(1,q,r)$ is the diagonal linear transformation with $1$, $q$ and $r$ on
the main diagonal, then
$$g^*(d)=\psa{k-1,1,0}1q\lambda+\psa{k,0,0}3\mu r\inv
+\psa{k-2,1,1}3qk\lambda,$$ so that if neither $\lambda$ nor $\mu$
vanishes  the codifferential is equivalent to
\begin{math}
d_k^\sharp=\psa{k-1,1,0}1+\psa{k,0,0}3+\psa{k-2,1,1}3k.
\end{math}
When $\lambda=0$, we obtain the codifferential
\begin{math}
d_k^*=\psa{k,0,0}3.
\end{math}
The case $\mu=0$ is simply a duplicate of the codifferential
$d_k(\lambda:k\lambda)$.

\subsection{Extensions of the Codifferentials}
We now study the deformations and extensions of each of these
types of codifferentials. We will show in sections
\ref{section5.7} and \ref{section5.8} that the cohomology $H^n_1$
vanishes for $n>k$ for the codifferentials $d_k^\star$ and
$d_k^\sharp$.  Because of this, they have no nonequivalent
extensions. In general, this will be true for the codifferentials
$d_k(\lambda:\mu)$, except for some special cases.

Let $d=d_k(\lambda:\mu)$. The matrix of $D:C^l_1\ra C^{k+l-1}_2$
is $\left[\begin{smallmatrix}0&\lambda l
-\mu&0\end{smallmatrix}\right]$, so unless $\mu=\lambda l$, the
dimension of the space $Z^l_1$ of $l$-cocycles is 2. Thus the
codifferential $d_k(1:l)$ is special when $l$ is any positive
integer, because the dimension of the space of cocycles is one
higher than usual.

The matrix of $D:C^{l}_0\ra C^{k+l-1}_1$ is
$\left[\begin{smallmatrix}
\lambda(k-l-1)&0&\lambda&0\\
0&0&0&0\\
\mu(k-2)&0&\mu&-\lambda(l-1)
\end{smallmatrix}\right]$. When $\lambda=0$ or $k=l$, the rank of this
matrix is equal to 1; otherwise the rank is 2. Thus generically,
when $l>k$ we have $H^n_1(d_k(\lambda:\mu))=0$ for all $n>k$. The
only cases where nontrivial extensions occur are for the
codifferentials $d_k(0:1)$ and $d_k(1:l)$, when $l>k$.

\subsubsection{Extensions of $d_k(0:1)$}
In this case we have
\begin{align*}
H^1_1=&\langle\psa{0,1,0}1\rangle\\
H^n_1=&\langle\psa{n-1,1,0}1,\psa{n-2,1,1}3\rangle,\quad 1<n<k\\
H^n_1=&\langle\psa{n-1,1,0}1\rangle,\quad n\ge k.
\end{align*}
We can extend $d_k(0:1)$ to the codifferential
$\psa{k-2,1,1}3+\psa{l-1,1,0}1x$. When $x\ne0$, by applying a
linear automorphism, we can transform this codifferential to one
of the same form with $x=1$. Thus, up to equivalence, for each
$l>k$, we obtain an extension
\begin{equation}
d_{k,l}=\psa{k-2,1,1}3+\psa{l-1,1,0}1.
\end{equation}
It can be shown that $H^n(d_{k,l})=0$ for $n>l$, unless $n=2l-k$.
Since the cohomology does not vanish completely, it is natural to
suspect that one can add a nontrivial term of degree $2l-k$. In
fact,  any extension of $d_{k,l}$ must be of the form
\begin{equation*}
d^e=d_{k,l}+\psa{n-1,1,0}1\alpha+\psa{n-2,1,1}3\beta+\ho,
\end{equation*}
since the third nonvanishing term in such an extension has to be a
cocycle with respect to $d_k$. Applying the automorphism
$\exp(-\phd{n-k,1,0}2\beta)$ to $d_e$ will replace the term
$\beta$ term in $d^e$ with higher order terms, so we can assume it
is zero. If $\eta=\phd{n-l+1,0,0}1-\phd{n-l,1,0}2(k-2) $, then
$D(\eta)=0$ and $[\psa{l-1,1,0}1,\eta]=-\psa{n-1,1,0}1(n+k-2l)$,
so as long as $n\ne 2l-k$, we can eliminate the $\alpha$ term in
$d^e$ as well. We cannot eliminate the term in degree $2l-k$,
which is exactly the degree in which the cohomology does not vanish.
Thus we obtain a nontrivial extension of $d_{k,l}$ given by
\begin{equation}
d_{k,l}(\alpha)=d_{k,l}+\psa{2l-k-1,1,0}1\alpha.
\end{equation}
Moreover,  a diagonal automorphism $g=\diag(p,q,t)$  leaves
$d_{k,l}$ invariant precisely when $p^{l-k}=1$, and in this case
$g$ leaves $\psa{2l-k,1,0}1$ invariant. From this it follows that
the extensions $d_{k,l}(\alpha)$ are nonequivalent for different
values of $\alpha$. Thus we obtain a family of nonequivalent
extensions of $d_k(0:1)$.

But this is the end of the story, because it can be shown that any
extension of $d_{k,l}(\alpha)$ is equivalent to it.  Note this is
not surprising, because $H^n_1(d_{k,l}(\alpha))=0$ for $n>2l-k$.

\subsubsection{Extensions of $d_k(1:l)$, where $l>k$}
Since $\psa{l,0,0}3$ is a nontrivial cohomology class, we obtain a
nontrivial extension of $d_k(1:l)$ by adding a multiple of this
cocycle, and by the usual argument, we can assume the coefficient
of the added term is one.  Thus we obtain the nontrivial extension
\begin{equation}
d^\sharp_{k,l}=\psa{k-1,1,0}1+\psa{l,0,0}3+\psa{k-2,1,1}3l.
\end{equation}
The matrix of $\ad_{\psa{l,0,0}3}:C^n_0\ra C^{l+n-1}_1$ is
$\left[\begin{smallmatrix}
0&1&0&0\\
l&0&0&-1\\
0&l&0&0
\end{smallmatrix}\right]$, and the matrix of
$\ad_{\psa{l,0,0}3}:C^n_1\ra C^{l+n-1}_2$ is
$\left[\begin{smallmatrix}l&0&-1\end{smallmatrix}\right]$. It
follows that $d_{k,l}^\sharp$-cocycles of degree larger than $l$
have the form $\psa{n-1,1,0}1+\psa{n-1,1,1}3l$, and these terms
are always $d_{k,l}^\sharp$-coboundaries. Thus
$H^n(d^\sharp_{k,l})=0$ vanishes for $m\ge l$, which suggests that
the extensions $d^\sharp_{k,l}$ are the only nontrivial extensions
of $d_k(1:l)$. In fact, if we write an extension of
$d_{k,l}^\sharp$ by adding terms of degree $n$ or higher, then it
has to be of the form
\begin{equation*}
(d_{k,l}^\sharp)^e=\d_{k,l}^\sharp+\psa{n-1,1,0}1\alpha+\psa{n-2,1,1}3\beta
+\ho,
\end{equation*}
since the degree $n$ term must be a $d_k(1:l)$-cocycle. But the
term of degree $n$ can be eliminated by applying an exponential,
so we don't obtain any higher order nonequivalent extensions. Note
that if $l\le k$, then  $d_k(1:l)$ has no nontrivial extensions.

\subsection{Deformations of the nonequivalent codifferentials}
As a result of our study of extensions, we can now give a complete
classification of all nonequivalent codifferentials on $W$. Let us
relabel $d_k^\sharp$ as $d_{k,k}^\sharp$, because it fits nicely
into the pattern $d_{k,l}^\sharp$, and we shall see later that it
deforms in a similar manner to elements of this form. The complete
list of nonequivalent codifferentials is given by
\begin{equation}\label{codiffs2}
\begin{array}{lcl}
d_k(\lambda:\mu)&=&\psa{k-1,1,0}1\lambda+\psa{k-2,1,1}3\mu\\
d_k^*&=&\psa{k,0,0}3\\
d_{k,l}(\alpha)&=&\psa{k-2,1,1}3+\psa{l-1,1,0}1+\psa{2l-k-1,1,0}1\alpha,
\quad k\le l\\
d^\sharp_{k,l}&=&\psa{k-1,1,0}1+\psa{k-2,1,1}3l
+\psa{l,0,0}3,\quad k\leq l.
\end{array}.
\end{equation}
We now construct miniversal deformations for each of these
codifferentials and
 then determine which codifferentials
arise as we substitute values of the parameters into the
miniversal deformation. This tells us how the codifferentials
deform locally, and gives us the gluing relations on the moduli
space.

\subsubsection{Deformations of $d_k(\lambda:\mu)$: Generic Case}
In this section, we suppose that $\lambda\ne0$ and
$\mu/\lambda\notin\Z_+$.  Then the cohomology is given by
\begin{align*}
H^1_1=&\langle\psa{0,1,0}1\rangle\\
H^n_1=&\langle\psa{n-1,1,0}1,\psa{n-2,1,1}3\rangle&1<n<k\\
H^k_1=&\langle\psa{k-1,1,0}1\rangle\\
H^n_1=&0,&n>k.
\end{align*}
As a consequence,  the universal infinitesimal deformation is
given by
\begin{align*}
d^\infty_k(\lambda:\mu)=&
d_k(\lambda:\mu)+\psa{m-1,1,0}1s_m+\psa{n-2,1,1}3t_n,1\le m<k,1<n\le k\\
&=d_k(\lambda:\mu+t_n)+\sum_{n=1}^{k-1}d_n(s_n:t_n).
\end{align*}
Furthermore, this deformation is miniversal since its self-bracket
vanishes. In other words, the miniversal deformation is given by the same formula
above, interpreting the parameters as lying in the formal power series algebra $\C[[s_n,t_n]]$.
Because the bracket vanishes, there are no relations on the base of the miniversal deformation.

If we substitute some values for the parameters in the formula above, we want
to determine the equivalence class of the resulting
codifferential. Let $M$ be the least value of $m$ so that
$s_m\ne0$, if it exists, and similarly, let $N$ be the least value
of $n$ so that $t_n$ does not vanish. If $s_m$ vanishes for all
$m$ and $N<k$, then the deformation is equivalent to
$\displaystyle d_{N,k}(-\tfrac{t_k+\mu}\lambda)$. When $s_m$ vanishes for all
$m$ and $N=k$, the deformation is equal to $d_k(\lambda:\mu+t_n)$,
which lies in the family.

If $M\le N$ or $t_n=0$ for all $n$, then the deformation is
equivalent to $d_M(s_M:t_M)$, since all the higher order terms are
coboundaries with respect to the leading term $d_M(s_M:t_M)$ in
$\dinfty$. Since the coordinates $(s_M:t_M)$ are projective, there
are deformations to $d_M(s_M:t_M)$ for arbitrarily small values of
the parameters, which means that the deformations to
$d_M(s_M:t_M)$ are jump deformations.

If $M>N$, then the deformation is equivalent to $d_{N,M}(\beta)$,
where $\beta$ depends rationally on the parameters $t_n$ and
$s_{m}$. The denominator is $s_M^{M-N+1}$, and the numerator is a
polynomial in $s_M$ of degree $M-N$, with coefficients given by
polynomials  in the parameters $t_n$ and $s_m$ for $m>M$, with no
constant term in any of the coefficients. For example, if $k=5$,
$N=3$ and $M=5$, then $\displaystyle\beta={\frac {-{s_{{5}}}^{2}t_{{5}}+
\left( \lambda\,t_{{3}}+s_{{6}}t_{{4}}
 \right) s_{{5}}-{s_{{6}}}^{2}t_{{3}}}{{s_{{5}}}^{3}}}$.
Therefore, every possible
value of  $\beta$ arises along curves with arbitrarily small
values of the parameters. As a consequence, there are jump
deformations to $d_{N,M}(\beta)$ for all $\beta$.

The only deformations which do not occur as jump deformations are
the deformations $d_k(\lambda:\mu+t_n)$ along the family. Thus the
pattern is that we deform along the family $d_k(\lambda:\mu)$ of
the same order of codifferential, and to lower order objects of the form
$\d_n(\alpha,\beta)$ for all $(\alpha:\beta)$, as well as to all
$d_{n,m}(\beta)$ for $n<k$ and $m<k$.

\subsection{Deformations of $d_k(0:1)$}
The miniversal deformation is
\begin{equation*}
d^\infty_k(0:1)= d_k(0:1)+\psa{m-1,1,0}1s_m+\psa{n-2,1,1}3t_n,1\le
m<\infty,1<n<k.
\end{equation*}
The deformation picture is the same as in the generic case, with
the following exceptions.  If $s_m=0$ for all $m$, then the
deformation is equivalent to $d_N(0:1)$. The new case of interest
is when $M>k$. Then the deformation is equivalent to $d_{N,M}(\beta)$,
where $\beta$ is a rational function with denominator $s_M^{M-N+1}$,
the numerator is a polynomial in $s_M$ of degree less than $M-N+1$,
with coefficients given by polynomials without constant terms in the parameters $t_n$ and
$s_m$ for $m>M$.  For example, if $k=7$, $N=6$, and $M=8$, then
$\displaystyle\beta={\frac {s_{{10}}t_{{6}}s_{{8}}+s_{{9}}t_{{7}}s_{{8}}+s_{{9}}\mu\,s_{{8
}}-{s_{{9}}}^{2}t_{{6}}}{{s_{{8}}}^{3}}}$.  As a consequence, we see that there
are jump deformations to $d_{N,M}(\beta)$ for all values of $\beta$.
Note that the
deformation along the family is given by $d_k(s_k:1)$ instead of
$d_k(\lambda:\mu+t_k)$, and this deformation is the only one which
is not attainable as a jump deformation.

\subsection{Deformations of $d_{k,l}(\alpha)$}
The miniversal deformation is given by
\begin{equation*}
d^\infty_{k,l}(\alpha)=
d_{k,l}(\alpha+s_{2l-k})+\sum_{n=1}^{l-1}d_n(s_n,t_n),
\end{equation*}
where $t_1=0$ and $t_n=0$ for $n>k$. The picture of the deformations is similar to
that of $d_k(0:1)$.
When $M=k$ and $N=k$ or $N$ does not exist, we get the deformation $d_k(s_k:t_k+1)$.
Otherwise,
when $M< k$ and $M\le N$ or $N$ does not exist, the deformation is equivalent to $d_M(s_M:t_M)$.
Also, as we explain below,
deformations of the form $d_{N,M}(\beta)$ arise only when $N<M<l$.

When $N<M\le k$, or when $N\le k$ and $M=k+1$
we get jump deformations to $d_{N,M}(\beta)$
for all values of $\beta$. We also get
jump deformations to $d_{k,k+1}(\beta)$ for all $\beta$, which corresponds to
the case when $N$ does not exist, and $M=k+1$.
If $s_m=0$ for all $m<2l-k$ then the deformation is equivalent to
$d_{N,l}((1+t_N)(\alpha+s_{2l-k}))$ if $N<k$, and  to $d_{k,l}(\alpha+s_{2l-k})$, when $N$ does
not exist.
These are smooth
deformations to a neighborhood of $d_{N,l}(\alpha)$ (or $d_{k,l}(\alpha)$), not jump deformations.

When $N<k\le M$ and $M\neq k+1$, or $N=k$ and $M>k+1$, the deformation is again of the form $d_{N,M}(\beta)$,
where again, $\beta$ is a rational function of the parameters with denominator
$s_M^{M-N+1}$ and the numerator is a polynomial in $s_M$ of degree $M-N$ or less, but this time,
one of the coefficients of the polynomial is simply $t_k+1$. The case where $N$ does not exist
is similar, but this time it is to $d_{k,M}(\beta)$ and one of the coefficients of  the polynomial
is  1.  This means that
the numerator is bounded away from zero for small values of the parameters, while the
denominator is close to zero. For example, if $k=6$, $l=9$, $N=4$ and $M=8$, then
the deformation is equivalent to
$\displaystyle{d_{4,8}\left({\frac {(\alpha+s_{12})\,t_{{4}}{s_{{8}}}^{3}- \left( t_{{6}}+1 \right) {s_{
{8}}}^{2}+s_{{8}}t_{{5}}-t_{{4}}}{{s_{{8}}}^{5}}}
\right)}$. As a consequence, we obtain smooth deformations to
$d_{N,M}(\beta)$ in a neighborhood of $\beta=\infty$.

What is the meaning of the deformations which occur only for large
values of $\beta$? If $l>k+1$, then we can think of the  $d_{k,l}(\alpha)$ as points lying
in an affine plane in $\P^1(\C)$, parameterized by the coordinate $\alpha$,
with $d_{k,l+1}(\beta)$ glued in as the point at
infinity.  Furthermore, each of these families $d_{k,l}(\alpha)$ is
glued in as the point at infinity for each smaller one
(\ie, for each $d_{m,n}(\beta)$ with $m<k\leq n<l$).
Thus we can
always deform to lower order objects on the sequence $d_{k,l}(\alpha)$
to points near infinity, that is, deformations with large values
of $\beta$. The exception is that $d_{k,k+1}(\alpha)$ has jump deformations
from $d_{k,l}(\beta)$ for all values of $\alpha$ and $\beta$ and all $l>k+1$.

\subsection{Deformations of $d_k(1:l)$}
The universal infinitesimal deformation is
\begin{align*}
\d_k^{\text{inf}}(1:l)=
d_{k}(1:l+t_k)+\sum_{n=1}^{k-1}d_n(s_n,t_n)+\psa{l,0,0}3r,
\end{align*}
where $t_1=0$. This time, the self bracket does not vanish. In
fact,
\begin{equation*}
\tfrac12[\d_k^{\text{inf}}(1:l),\d_k^{\text{inf}}(1:l)]=
\sum_{n=1}^kr(ls_n-t_n).
\end{equation*}
Some of these terms are coboundaries, so we can add some preimages
of coboundaries to cancel them. The terms which are not
coboundaries will lead to the relations on the base of the
miniversal deformation. We will need to add some terms to the
infinitesimal deformation to get the miniversal deformation.
These terms can be determined recursively, as follows.  Suppose
that
\begin{equation*}
d^\infty_k(1:l)=\d_k^{\text{inf}}(1:l)+\sum_{m=1}^{l-1}\psa{m,0,0}3x_m,
\end{equation*}
in other words, it is obtained by adding terms involving the
preboundaries $\psa{m,0,0}3$. Then we have
\begin{align*}
\tfrac12[d^\infty,d^\infty]= \sum_{m=1}^{l-1}\sum_{n=1}^k&
\phd{k+m-2,1,0}3x_m(m-l)+
\phd{l+n-2,1,0}3r(ls_n-t_n)\\&+\phd{n+m-2,1,0}3x_m(ms_n-t_n).
\end{align*}
These equations can be solved for the coefficients $x_m$, which
will be given by rational expressions in the parameters $r$, $s_n$
and $t_n$. In fact $x_{l-1}=\frac{r(ls_{k-1}-t_{k-1})}{1+t_k}$,
and the others can be found by a downward recursion.  Note that as
a formal power series in the parameters, $x_{l-1}$ is of order 2;
the other $x$'s will also be of at least this order.

There are $k+l-1$ terms in the self-bracket of $d^\infty$, of
which the terms of degrees $1$ through $k-1$ and degree $k+l-1$
correspond to nontrivial cocycles. The other $l-1$ terms are used
to solve for the $x_m$ recursively.

As a simple example, let $l=2$ and $k=3$. Then we need to add only
$x_1=\frac{r(2s_2-t_2)}{1+t_3}$. There are 3 relations on the
base,
\begin{equation*}
rt_3=0,\qquad 2rs_1+x_1(s_2-t_2)=0,\qquad x_1s_1=0.
\end{equation*}
These relations must be satisfied by any codifferential which is
obtained from the miniversal deformation by substituting values in
for the parameters, and any substitution for the parameters
satisfying these relations will determine an actual
codifferential.  To understand the information encoded in the
miniversal deformation, we need to analyze the equivalence classes
of the resulting codifferentials we obtain in this manner.  For
our example, we see that a deformation with $x_1\ne0$ will be
equivalent to $d_1^*$. If $x_1=0$, and $r=0$, the deformation is
equivalent to either $d_1(1:0)$, $d_2(s_2:t_2)$ or $d_3(1:l+t_3)$.
On the other hand, if $r\ne0$, then $s_1=t_3=0$ and $2s_2=t_2$.
This give rise to $d_{2,2}^\sharp$.

The general pattern is that $d_n^\star$ will arise as a
deformation of $d_k(1:l)$; in addition, $d_{n,m}^\sharp$ will
appear when $n\le k$ and $m\le l$.

\subsection{Deformations of $d_{k,l}^\sharp$}\label{section5.7}
The cohomology is given by
\begin{align*}
H^1_1=&\langle\zeta=\psa{0,1,0}1(k-1)+\psa{1,0,0}3l\rangle\\
H^n_1=&\langle d_n(1:l),d_{n,n+l-k}\rangle, 1<n<k\\
H^n_1=&0, n\ge k
\end{align*}
The universal infinitesimal deformation is given by
\begin{align*}
(d_{k,l}^\sharp)^\text{inf}=&d_{k,l}^\sharp+\zeta r +d_m(1:k)s^m
+d_{n,n+l-k} t^n,\\& 1<m<k, 1<n<k
\end{align*}
The term $\zeta$ coming from $H^1$ does not fit the nice pattern
we have been observing, and indeed, when we compute the
self-bracket of $(d_{k,l}^\sharp)^\text{inf}$, the lowest degree
term is multiplied by
 $2lr^2(k-1)(l-k+1)$, which means that in any actual
deformation, we must have $r=0$. In classical language, this means
that $d_{k,l}^\sharp+\zeta r$ does not extend to a second order
deformation.  Let us assume that $r=0$ and then compute the self
bracket. We obtain
\begin{align*}
\tfrac12[(d_{k,l}^\sharp)^\text{inf},(d_{k,l}^\sharp)^\text{inf}]=
\phd{n-2+l-k+m,1,0}3s^mt^n(n-k).
\end{align*}
It is easy to see that $d_{k,l}^\sharp$ deforms to the lower order
codifferentials $d_{n,n+l-k}^\sharp$ as well as the codifferentials
$d_{n,l}^\sharp$, when $n=2\dots k-1$. What is less obvious is that
we also may obtain codifferentials of the form $d_n^*$ for $n=1\dots
k-1$ and codifferentials of the form $d_{n,m}^\sharp$ for $1<n\le
m<k$.

In general, it is complicated to describe precisely which lower order
codifferentials will arise. For example, if $k=3$ and $l=4$,  the
codifferentials $d_1^*$ and $d_{2,2}^\sharp$ occur, but not
$d_2^*$. Thus, all deformations are to codifferentials of the type
listed above, but not all of the possibilities actually occur for
a particular codifferential.

\subsection{Deformations of $d_k^\star$}\label{section5.8}
The cohomology is somewhat different for this codifferential. We
have
\begin{align*}
H^1_1=&\langle\psa{1,0,0}3\rangle\\
H^n_1=&\langle\psa{n,0,0}3,\psa{n-1,1,0}1+\psa{n-2,0,1}3n\rangle& 1<n<k\\
H^k_1=&\langle\psa{k-1,1,0}1+\psa{k-2,0,1}3k\rangle\\
H^n_1=&0,&n>k.
\end{align*}
As a consequence, we see that there are no extensions of
$\d_k^\star$, and that it deforms into the codifferential
$d_{k,k}^\sharp$, of the same order. The universal infinitesimal deformation is given
by
\begin{equation*}
(d_k^*)^{\text{inf}}=
d_k^*+\psa{m,0,0}3s_m+(\psa{n-1,1,0}3+\psa{n-2,1,1}3k)t_n, 1\le
m<k,1<n\le k.
\end{equation*}
We compute the self bracket
\begin{equation}
\tfrac12[(d_k^*)^{\text{inf}},(d_k^*)^{\text{inf}}]=\phd{m+n-2,1,0}3s_mt_n(m-k).
\end{equation}
Since coboundary terms occur in the bracket, we have to add some
terms of the form $\psa{r-1,1,0}1x_r$, in fact, exactly $k-1$ of
them to correspond to the $k-1$ coboundary terms in the relations
above. Note that in the end, there will be $k-1$ relations on the
base of the miniversal deformation.  In this case, we do not need to
solve all of the equations in order to determine what
codifferentials arise as deformations. It is clear that we can
deform to $d_m^\star$ for $1\le m<k$ as well as to $d_{n,k}^\sharp$
for $1<n\le k$.  But it may occur that the lowest $m$ term added is
of larger degree than the lowest $n$ term.  In this case, we would
expect that the deformation is equivalent to $d_{n,m}^\sharp$.
However, this would not be true without the addition of the $x_r$
terms, which fortunately are added in the versal deformation. These
terms are precisely what is necessary in order for the deformed
codifferential to be equivalent to $d_{n,m}^\sharp$. Thus, we can
also deform to the codifferentials $d_{n,m}^\sharp$ for  $1<n\le m< k$.
\begin{figure}
\begin{center}
\epsfxsize=4.9in \epsfysize=4.6in \epsfbox{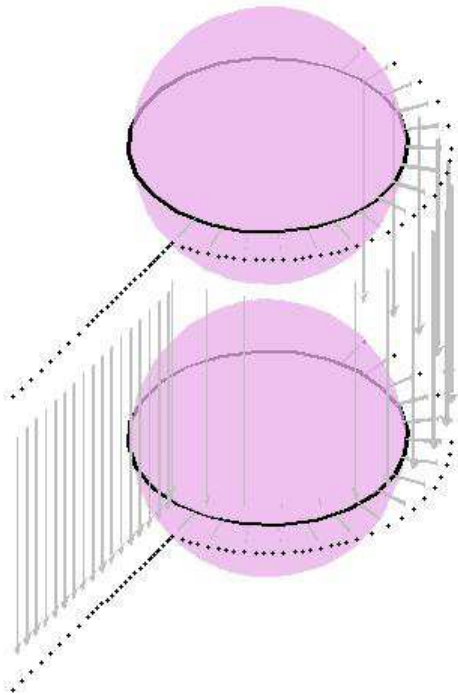}
\caption{{\bf The moduli space of codifferentials on } 
$\mathbf{\W_{-1}\oplus W_{0}\oplus W_{1}}$.
This picture is explained in section 
\ref{long_caption}.} \label{zeroonegone_label}
\end{center}
\end{figure}

\subsection{Description of Figure \ref{zeroonegone_label}}\label{long_caption}

The upper large sphere represents $d_k(\lambda:\mu)$, and the large sphere below
it is $d_{k-1}(\lambda:\mu)$.
Although the diagram only shows two distinct strata (at levels $k$ and $k-1$), there
is actually an infinitely long chain of them.
We will now describe the $k^{th}$ stratum
(shown in the upper half of the picture) in detail.

The large sphere, $d_k(\lambda:\mu)$, can be thought of as the
complex projective space $\P^1(\C)$.
The point on the equator which has a long spike extending above it is $d_k(0:1)$,
the point at infinity.
The equator also includes all points $d_k(1:l)$.
The points above the equator (except for the long spike) 
represent $d^{\sharp}_{k,l}$ (for each $l\in\mathbb Z$).
Recall that when $l\in\mathbb Z$, 
$d_k(1:l)$ has a jump deformation to $d^{\sharp}_{k,l}$,
which is represented by the arrows radiating out from the equator.
Note that there are also many other deformations of this sort which are not 
explicitly shown in our diagram (for the sake of simplicity).
Specifically, 
$d_k(1:l)$ has jump deformations to 
$d^{\sharp}_{m,n}$ for all $m\leq k$ and $n\leq l$.

The long spike above the point at infinity of the
large sphere represents the family $d_{k,l}(\alpha)$.  
The \lq\lq smallest\rq\rq\ element of the family on level $k$
is $d_{k,k+1}(\alpha)$, which is on the end of the spike which
is farthest from the large sphere $d_k(\lambda:\mu)$.
The entire spike consists of an infinite sequence,
$\{d_{k,l}(\alpha)|l>k\}$, which gets closer to the
sphere as $l$ increases.

Since each member of the family $\{d_{k,l}(\alpha)|l>k\}$ 
is parameterized
by a complex number $\alpha$ and has a point at 
infinity, we can visualize each $d_{k,l}(\alpha)$ 
as a sphere.
Specifically, we can think of $d_{k,l+1}(0)$
as the point at infinity on $d_{k,l}(\alpha)$,
since $d_{k,l+1}(0)$ deforms smoothly to 
$d_{k,l}(\beta)$ in a neighborhood of $\beta=\infty$.
From this point of view, one could visualize 
$\{d_{k,l}(\alpha)|l>k\}$ as an infinitely long string of pearls,
as indicated by the diagram.
However, this pretty picture is an oversimplification,
since we can actually deform smoothly 
to $d_{k,l}(\beta)$ in a neighborhood of $\beta=\infty$
from $d_{m,n}(\alpha)$ for all $k\leq m\leq l<n$,
and for all $\alpha$.

From each $d_{k,l}(\alpha)$,
there is a smooth deformation to a neighborhood of
each $d_{n,l}(\alpha)$ \lq\lq below\rq\rq\ it
(\ie, for each $n<k$), represented by the
downward arrows descending from the spike.
In the interest of keeping the diagram
uncluttered, there are many other deformations
which are not explicitly represented by arrows,
but which are described in the next paragraph.

From each point on  $d_k(\lambda:\mu)$,
there is a jump deformation to each \lq\lq lower\rq\rq\
$d_{m,n}(\beta)$ (in other words, for each $m<n\leq k$).
Additionally, there is a jump deformation from 
$d_{k,l}(\alpha)$ to each $d_{m}(\alpha:\beta)$ 
which is at the same level or lower, \ie, for each $m\leq k$.
From each point $d_k(\lambda:\mu)$ on the large sphere,
there exist smooth deformations to every other point 
$d_k(\alpha:\beta)$ on the same sphere.
There also exist jump deformations to all points
on the \lq\lq lower\rq\rq\ spheres.  In other words, 
there is a jump deformation from $d_k(\lambda:\mu)$ 
to $d_n(\alpha:\beta)$ for all $n<k$.

For the sake of simplicity, $d^*_k$ is omitted from the diagram,
although it could be viewed as an extra point at level $k$.
From $d^*_k$, there exist jump deformations to $d^\sharp_{k,k}$, 
to $d^*_m$ (for all $1\leq m\leq k$), and to 
to $d^\sharp_{m,n}$ (for all $1\leq m< n< k$).
Furthermore, there are also jump deformations to 
$d^*_k$ from each $d^\sharp_{m,n}$ for $k\leq m$
which are not shown here.

\subsection{Comparison to \zt-graded codifferentials}
Given the complexity of the description of the moduli space for
the \Z-graded codifferentials, it would not be surprising if the
\zt-graded picture were different.  However, as we shall now show,
the deformation picture in the \zt-graded case is identical to the
\Z-graded picture.

If we consider more general codifferentials, then we have to take
into account the possibility that there are terms coming from
$C^k_{-1}$ in a degree $k$ codifferential. Thus we can express a
general odd cochain in the form
\begin{align*}
d=\psa{k-1,1,0}1a_1+\psa{k,0,0}3b_1+\psa{k-2,1,1}3c_1+
\psa{k-1,0,1}1a_2+\psa{k,0,0}2b_2+\psa{k-2,1,1}3c_2.
\end{align*}
It is easy to check that the condition $[d,d]=0$ gives, among
other relations, the requirements that $b_2a_1+a_2b_1=0$ and
$c_2a_1+a_2c_1=0$. If we consider a linear automorphism given by
the matrix $ \left[
\begin{smallmatrix}
1&0&0\\
0&r&s\\
0&t&u
\end{smallmatrix}
\right] $, then $g^*(d)$ is of the form

\begin{equation*}
g^*(d)=\psa{k-1,1,0}1a+\psa{k,0,0}3b+\psa{k-2,1,1}3c
\end{equation*}
precisely when
\begin{align*}
-sb_1+ub_2=0,\quad&sa_1+ua_2=0,\quad-sc_1+uc_2=0.
\end{align*}
Such an automorphism always exists, so that any codifferential of
degree $k$ is equivalent to a codifferential lying in $C^k_1$. In
fact, if we consider the standard basis
$\{e_1^k,e_1^{k-1}e_2,e_1^{k-1}e_3,e_1^{k-2}e_2e_3\}$ of $S^k(W)$,
then the matrix of $d$ with respect to this basis is
$A=\left[ \begin {smallmatrix} 0&{\it a_1}&{\it a_2}&0\\
\noalign{\medskip}{\it b_2}&0&0&{\it c_2}\\
\noalign{\medskip}{\it b_1}&0&0&{\it c_1}
\end {smallmatrix} \right]$, while the matrix of $g:S^k(W)\ra S^k(W)$ with
respect to the standard basis is $Q= \left[
\begin {smallmatrix} {q}^{k}&0&0&0\\
\noalign{\medskip}0&r{q}^{k-1}&s{q}^{k-1}&0\\
\noalign{\medskip}0&t{q}^{k-1}&u{q}^{k-1}&0\\
\noalign{\medskip}0&0&0& \left( ur-st \right) {q}^{k-2}
\end {smallmatrix}
\right] $, so that the matrix of $g^*(d)$ is
$$A'=G\inv A Q=
\left[
\begin {smallmatrix} 0&{q}^{k-2} \left( {\it a_1}\,r+t{\it a_2} \right)
&{q}^{k-2} \left( s{\it a_1}+{\it a_2}\,u \right) &0\\
\noalign{\medskip} {\frac { \left( u{\it b_2}-s{\it b_1} \right)
{q}^{k}}{ur-st}} &0&0&{q}^{k-2} \left( u{\it c_2}-s{\it c_1}
\right)
\\\noalign{\medskip}
{\frac { \left( -t{\it b_2}+r{\it b_1} \right) {q}^{
k}}{ur-st}}&0&0&{q}^{k-2} \left( -t{\it c_2}+r{\it c_1} \right)
\end {smallmatrix} \right]
.$$ Note that this matrix represents an element in $C^k_1$
precisely when the second row and third column vanish.  Now
suppose that $d$ already lies in $C^k_1$.  It follows that
$g^*(d)$ lies in $C^k_1$ exactly when $s=0$.  Thus when studying
even automorphisms preserving $C^k_1$, we don't get a restriction
just to diagonal automorphisms, but we do have a simple
characterization of the automorphisms.  Moreover, the action of
$g$ on $d$ depends only on the diagonal components, which are the
degree zero components, so the action is the same as the degree
zero automorphisms acting in the \Z-graded case.  As a
consequence, we obtain the same equivalence classes of \zt-graded
codifferentials as \Z-graded codifferentials.

To understand the moduli space of \zt-graded codifferentials, we
need to study the deformations.  We already understand the
coboundaries of degree 1 and degree zero cochains,  but we need to
study the action of the coboundary operators on degree $-1$ and
degree $-2$ (the action on the degree 2 part is trivial). Let us
examine the situation for each of our codifferentials.

\subsubsection{The codifferential $d_k(\lambda:\mu)$ and $d_{k,l}$}
$D:C^1_{-1}\ra C^{k+l-1}_0$ is given by $
\left[\begin{smallmatrix}
0&\lambda&0\\
\lambda(l-k)+\mu&0&\lambda\\
0&l\lambda&0\\
0&\mu&0
\end{smallmatrix}
\right] $, while that of $D:C^{l}_{-2}\ra C^{k+l-1}_{-1}$ is
$\left[\begin{smallmatrix} \lambda\\0\\-(\lambda(l-1)+\mu)
\end{smallmatrix}\right]$.
For $l>1$, there is exactly one cocycle of degree -1,
$$\delta_l=\psa{l-1,0,1}1\lambda-\psa{l-2,1,1}2(\lambda((l-k)+\mu),$$
which is a coboundary when $l\ge k$.

Thus there are no additional extensions arising from the new
cocycles, but we do have to add some additional terms to the
miniversal deformation. In fact, we can generically express the
miniversal deformation in the form
\begin{equation*}
d_k^\infty(\lambda:\mu)=d_k(\lambda:\mu)+d_m(s^m,t^m)+\delta_lr^l,
1\le m<k,1<l<k.
\end{equation*}
It is easy to analyze that terms in the self bracket of
$d_k^\infty(\lambda:\mu)$ can only arise from the brackets of the
$d_m(s^m,t^m)$ terms with the $\delta_lr^l$ terms. We obtain
relations on the base of the miniversal deformation of the form
\begin{equation*}
\sum_{l+m=n+1}r^l(s^m(\lambda(m-k)+\mu)-t^n\lambda)=0,\quad 2\le
n\le2k-2.
\end{equation*}
What is the meaning of these relations? Let us use them to
determine the deformations.  First, let $l_0$ be the least $l$
such that $r^l\ne0$ and $m_0$ be the least $m$ such that either
$s^m$ or $t^m$ does not vanish. If $r^l$ vanishes for all $l$,
then we are already in the same situation as the \Z-graded case,
so we may assume otherwise. The first relation is just
$r^{l_0}(s^{m_0}(\lambda(m_0-k)+\mu)-t^{n_0}\lambda)=0$. If
$l_0<m_0$, then by applying an automorphism to
$d_k^\infty(\lambda:\mu)$, we transform it to one in which
$l_0>m_0$. If $m_0<l_0$, then the first relation says that
$\delta_{l_0}$ is a $d_{m_0}(s^{m_0},t^{m_0})$ cocycle, thus a
coboundary with respect to $\delta_{l_0}$, which means we can
eliminate it by applying an appropriate formal automorphism.
Finally, if $l_0=m_0$, we have exactly the relation between the
coefficients which makes it possible to transform the $m_0$ term
to one which lies entirely in $C^{m_0}_1$.  As a consequence, even
though the miniversal deformation in the \zt-graded case contains
additional terms, we obtain no new deformations in the \zt-graded
case.

The deformation picture is similar for the extended
codifferentials $d_{k,l}$.  We omit the detailed calculations, but
note that again the \zt-graded deformations are the same as for
the \Z-graded case.

\subsubsection{The codifferential $d_k^\star$}
$D:C^l_{-1}\ra C^{k+l-1}_0$ is given by$ \left[\begin{smallmatrix}
1&0&0\\
0&0&0\\
0&0&-1\\
k&0&0
\end{smallmatrix}
\right] $, while the matrix of $D:C^l_{-2}\ra C^{k+l-1}_{-1}$ is
$\left[\begin{smallmatrix}0\\1\\0\end{smallmatrix}\right]$. It
follows that the only cocycle is $\psa{l,0,0}2$, which is a
coboundary when $l\ge k$.  Thus we obtain no new extensions,  do
have to add terms to the miniversal deformation of $d_k^\star$.
\begin{align*}
(d_k^*)^\infty=&
d_k^*+\psa{1,0,0}2r_m+\psa{m,0,0}3s_m+(\psa{n-1,1,0}3+\psa{n-2,1,1}3k)t_n
\\&
1\le l,m<k,1<n\le k.
\end{align*}
We have
\begin{align*}
\tfrac12[(d_k^*)^\infty,(d_k^*)^\infty]=&\phd{m+n-2,1,0}3s_mt_n(m-k)
+\phd{l+n-1,0,0}1r_lt_n\\&+\phd{l+n-2,1,0}2lr_lt_n+\phd{l+n-2,0,1}3kr_lt_n.
\end{align*}
We already saw in studying the deformations of $d_k^\star$ in the
\Z-graded case that either $s_m=0$ for all $m$ or $t_n=0$ for all
$n$. In addition, we see that either $t_n=0$ for all $n$ or
$r_l=0$ for all $n$.  These last relations are exactly what we
need to be able to transform $(d_k^*)^\infty$, in the case when
$r_l\ne0$ for some $l$, to a codifferential lying entirely in
$C^\bullet_1$.  Thus we obtain no new deformations for $d_k^\star$
in the \zt-graded case.

\subsubsection{The codifferentials $d_k^\sharp$ and $d_{k,l}^\sharp$}
$D:C^l_{-1}\ra C^{k+l-1}_0$ is given by the matrix $
\left[\begin{smallmatrix}
1&1&0\\
l&0&1\\
0&l&-1\\
k&k&0
\end{smallmatrix}
\right] $,while the matrix of $D:C^l_{-2}\ra C^{k+l-1}_{-1}$ is
$\left[\begin{smallmatrix}1\\-1\\-(k+l-1)\end{smallmatrix}\right]$.
It follows that the only cocycles are of the form
$$\delta_l=\psa{l-1,0,1}1-\psa{l,0,0}2-\psa{l-2,1,1}2l,$$
which is a coboundary when $l\ge k$.  Thus we obtain no new
extensions, but as usual, do have to add terms to the miniversal
deformation of $d_k^\sharp$.
\begin{align*}
(d_k^\sharp)^\infty=&d_k^\sharp+(\psa{0,1,0}1(k-1)+\psa{1,0,0}3k)r
+d_m(1:k)s^m +d_n^\star t^n+\delta_lu^l,\\& 1<m<k, 1<n<k, 1<l<k.
\end{align*}
When the self bracket of $(d_k^\sharp)^\infty$ is computed, in
addition to the relations we saw in the \Z-graded case, we get
relations $\sum_{l+m=n+1}u^ls^m=0$, which means that either all
the $\delta_l$ or all the $d_m(1:k)$ terms must vanish in an
actual deformation.  But this is just what is necessary to be able
to transform a deformation with $\delta_l$ terms into one which
lies in $C^\bullet_1$. Moreover, no new deformations arise.  The
situation is similar for the codifferential $d^\sharp_{k,l}$.

As a consequence of these calculations, we have shown that the
\zt-graded moduli space and the \Z-graded moduli space for a
vector space of type $(-1,0,1)$ are completely identical,
including the deformations which glue the moduli space together.
As a result, Figure \ref{zeroonegone_label}, which illustrates how
the moduli space is glued together, applies to either the
\zt-graded or \Z-graded moduli spaces.

\section{Codifferentials on a $2|1$-dimensional space}
\label{sec2bar1} If $W$ is of the form $W=W_{2x}\oplus
W_{2x+1}\oplus W_{2x+2}$, then as a \Z-graded space it has type
$(2x,2x+1,2x+2)$, but its dimension is always $2|1$ as a
\zt-graded space. In order to list the even basis elements first,
we set $|e_1|=2x$, $|e_2|=2x+2$ and $|e_3|=2x+1$.

\begin{align*}
C^r_{2s}=&\langle
\phd{(r-1)(x+1)+s+1,-(r-1)x-s,0}1,\phd{(r-1)(x+1)+s,-(r-1)x-s+1,0}2,\\&
\phd{(r-1)(x+1)+s,-(r-1)x-s,1}3 \rangle\\
C^r_{2s+1}=&\langle
\psa{(r-1)(x+1)+s+1,-(r-1)x-s-1,1}1,\psa{(r-1)(x+1)+s,-(r-1)x-s,1}2,\\&
\psa{(r-1)(x+1)+s+1,-(r-1)x-s,0}3 \rangle
\end{align*}
In particular, and of most interest to us are the following.
\begin{align*}
C^r_0\!\!\!=&\langle
\phd{r+(r-1)x,-x(r-1),0}1,\phd{(r-1)(x+1),1-x(r-1),0}2,
\phd{(r-1)(x+1),-x(r-1),1}3 \rangle\\
C^r_1\!\!\!=&\langle
\psa{(r-1)(x+1)-1,-x(r-1)-1,1}1,\psa{(r-1)(x+1),-x(r-1),1}2,
\psa{(r-1)(x+1)+1,-x(r-1),0}3 \rangle
\end{align*}
There is a natural basis of $S^k(W)$, given by
$$\{e_1^k,e_1^{k-1}e_2,\mcom,e_2^k,e_1^{k-1}e_3,e_1^{k-2}e_2e_3,\mcom
e_2^{k-1}e_3\},$$ from which it follows that $S^n(W)$ has
\zt-graded dimension $(k+1)|k$.  If $A$ is the matrix of an odd
codifferential of degree $k$, then $A$ can be expressed in block
form as
$A=\left[\begin{smallmatrix}0&A_1\\A_2&0\end{smallmatrix}\right]$.
It can be checked that $A$ is a codifferential precisely when
either $A_1$ or $A_2$ vanishes.  A codifferential for which $A_2$
vanishes is said to be of the \emph{first kind}, and one for which
$A_1$ vanishes is of the \emph{second kind}.  More generally, a
codifferential with terms in more than one degree will still have
only terms of the first kind or second kind (see \cite{bfp1}). Thus the moduli
space of all codifferentials splits into disjoint subspaces, the
moduli spaces of codifferentials of the first and second kinds. We
will study these spaces separately.

\subsection{Codifferentials on a vector space of type (0,1,2)}\label{InjectiveExample}
Suppose $|e_1|=0$, $|e_2|=2$, and $|e_3|=1$. Then
$\phd{k,l,m}n=|e_n|-2l-m$ so
\begin{align*}
C^r_0=&\langle \phd{r,0,0}1,\phd{r-1,1,0}2,\phd{r-1,0,1}3 \rangle\\
C^r_1=&\langle \psa{r-1,0,1}2,\psa{r,0,0}3 \rangle.
\end{align*}
Up to equivalence there are only two codifferentials of degree
$k$, a codifferential of the first kind given by
$d=\psa{k-1,0,1}2$ and a codifferential of the second kind,
$d=\psa{k,0,0}3$.

\subsubsection{Codifferentials of the first kind}
Consider first the case
\begin{equation}
d_k=\psa{k-1,0,1}2.
\end{equation}
The cohomology of this codifferential is given by
\begin{align*}
H^n_1=&\langle \psa{n-1,0,1}2\rangle, \quad n<k\\
H^n_1=&0, \quad \text{otherwise.}
\end{align*}
Since $H^n_1$ vanishes for $n\ge k$ it follows that there are no
nontrivial extensions of $d$. Moreover
\begin{equation*}
\d_k^\infty=\psa{k-1,0,1}2+\sum^{k-1}_{m=1}\psa{m-1,0,1}2t_m,
\end{equation*}
so this means that the codifferentials of the first kind form a
discrete family with codifferentials of degree $k$ deforming into
codifferentials of degree $m$ for $m<k$ by jump deformations.

\subsubsection{Codifferentials of the second kind}
\begin{figure}[h,t]
\begin{center}
\epsfxsize=4.9in \epsfysize=3.0in \epsfbox{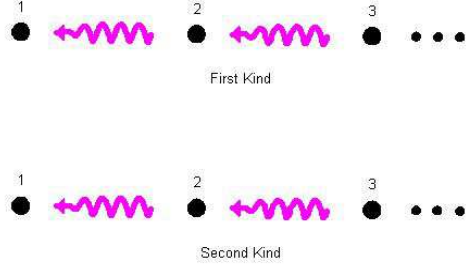}
\caption{{\bf The moduli space of codifferentials on }
$\mathbf{\W_0\oplus W_1\oplus W_2}$.
This moduli space contains two discrete families of
codifferentials.  The upper diagram shows the codifferentials
of the first kind, in which each $d_k$ has a jump
deformation to each \lq\lq lower\rq\rq\ $d_m$ (for all $m<k$).
The lower diagram represents the codifferentials of the 
second kind, which follow an analogous pattern.
} \label{caseone_label}
\end{center}
\end{figure}

For the codifferential of the second kind of degree $k$,
\begin{equation}d_k=\psa{k,0,0}3,\end{equation}
the pattern is similar. In this instance we have
\begin{align*}
H^n_1=&\langle \psa{n,0,0}3\rangle, \quad n<k\\
H^n_1=&0, \quad \text{otherwise}
\end{align*}
and
\begin{equation*}
\d_k^\infty=\psa{k,0,0}3+\sum^{k-1}_{m=1}\psa{m,0,0}3t_m.
\end{equation*}
Thus, the codifferentials of the second kind form another discrete
family. Figure \ref{caseone_label} is an illustration of the
moduli space of codifferentials.

Since the space of \zt-graded codifferentials is the same for this
example and for the example in the next section, we will defer the
analysis of the moduli space of \zt-graded codifferentials until
then.  However, it is important to note that for the \Z-graded
moduli space of type (0,1,2) considered here, the map to the
\zt-graded moduli space is injective, but not surjective.

\subsection{Codifferentials on a vector space of type (-2,-1,0)}\label{BadExample}
Suppose $|e_1|=-2$, $|e_2|=0$ and $|e_3|=-1$. Then
\begin{align*}
C^r_0=&\langle\phd{1,r-1,0}1,\phd{0,r,0}2,\phd{0,r-1,1}3\rangle\\
C^r_1=&\langle\psa{1,r-2,1}1,\psa{0,r-1,1}2,\psa{1,r-1,0}3\rangle,
\end{align*}
so we have a more complicated situation among the codifferentials
of the first kind, which are of the form
\begin{equation}
d_k(\lambda:\mu)=\psa{1,k-2,1}1\lambda+\psa{0,k-1,1}2\mu,
\end{equation}
where, for $k=1$ we omit the first term ($\lambda=0$). Using the
diagonal linear automorphism $g=\diag(p,q,r)$, it is easy to see
that
$$g^*(d_k(\lambda:\mu))=q^{k-2}rd_k(\lambda:\mu),$$ so the only
equivalent \Z-graded codifferentials are multiples. The situation
is quite different if we consider the even automorphism $g$ of $W$
which interchanges $e_1$ and $e_2$. This automorphism is even, but
not of degree 0 and
$$g^*(d_2(\lambda:\mu))=d_2(\mu:\lambda).$$ This means that
$d_2(\lambda:\mu)\sim d_2(\mu:\lambda)$ as \zt-graded
codifferentials, but not as \Z-graded codifferentials.  This
example shows that the map from the moduli space of \Z-graded
codifferentials to the moduli space of \zt-graded codifferentials
is not injective.  Note that the map from the \Z-graded moduli
space in the previous section is injective, because the \Z-graded
space in that case does not include codifferentials on which the
map fails to be injective.

The \Z-graded codifferentials of the first kind are labeled
projectively by points $(\lambda:\mu)\in\P^1(\C)$. The matrix of
$D:C^l_0\ra C^{k+l-1}_1$ is
\begin{equation*}
\left[
\begin{matrix}
\lambda&(k-2)\lambda&\mu(l-1)\\
\mu&\mu(k-l-1)&0\\
0&0&0
\end{matrix}
\right].
\end{equation*}
The rank of this matrix is equal to 2 except when $l=1$ or
$\mu=0$, when the rank is 1. This means, that the codifferential
$d_k(1:0)$ is special. We will study the extensions of the generic
case first.

\subsection{Extensions of $d_k(\lambda:\mu)$, $\mu\ne0$}
We have
\begin{align*}
H^1_1=&\langle\psa{0,0,1}2\rangle\\
H^n_1=&\langle\psa{1,n-2,1}1,\psa{0,n-1,1}2\rangle,\quad 1<n<k\\
H^k_1=&\langle\psa{0,k-1,1}2\rangle\\
H^n_1=&0,\quad n>k.
\end{align*}
Because of this, there are no nontrivial extensions of
$d_k(\lambda:\mu)$. Note that when $k=1$, we must have
$\lambda=0$, and then $H^n_1=0$ for all $n\ge 1$, so that
$d_1(0:1)$ has no extensions or deformations.

\subsection{Extensions of $d_k(1:0)$}
We have
\begin{align*}
H^1_1=&\langle\psa{0,0,1}2\rangle\\
H^n_1=&\langle\psa{1,n-2,1}1,\psa{0,n-1,1}2\rangle,\quad 1<n<k\\
H^n_1=&\langle\psa{0,n-1,1}2\rangle,\quad n\ge k.
\end{align*}

Because the cohomology does not vanish identically for $n>k$,
there are nontrivial extensions of $d(1:0)$.  If $l>k$, the
extension
\begin{equation*}
d_{k,l}=\psa{1,k-2,1}1+\psa{1,l-1,1}2x
\end{equation*}
is nontrivial when $x\ne 0$. Note that if $g=\diag(1,q,q^{2-k})$,
then
$$g^*(d_{k,e})=\psa{1,k-2,1}1+\psa{1,l-1,1}2q^{l-k}x,$$
 so that when $x\ne0$
we can assume that $x=1$.  Moreover, one can show that
$H^n_1(d_{k,l})=0$ for $n>l$, so that $d_{k,l}$ itself has no
nontrivial extensions.

\subsubsection{Deformations of co\-differentials of the first kind}
A codifferential of the first kind of degree $k$ is equivalent to
 $d_k(\lambda:\mu)$ or $d_{k,l}$.
 We will study the deformations of this space
in order to obtain a more detailed picture of the moduli space of
codifferentials of the first kind. If $k>1$, then
$H^k_1=\langle\psa{0,k-1,1}2\rangle$, so we always have a
deformation along the codifferentials of degree $k$; in other
words, we can deform along the family.  If $\mu\ne 0$, this
completes the picture in terms of deformations to codifferentials of the same order.
However, for $d_k(0:1)$, we can deform to $d_{k,l}$ for any $l>k$.
On the other hand,  $d_{k,l}$ deforms to $d_{k,m}$ for every
$k<m<l$, because $\psa{1,m-2,1}1$ remains a cohomology class for
$d_{k,l}$.

We can also deform to lower order codifferentials.  The codifferential
$d_k(\lambda:\mu)$ deforms to $d_m(\alpha:\beta)$ for any
$(\alpha:\beta)$ and all $1<m<k$, as well as to $d_{m,n}$ if
$1<m<n\le k$.  Similarly, $d_{k,l}$ deforms to $d_m(\alpha:\beta)$
for any $(\alpha:\beta)$ and all $1<m<k$, as well as to $d_{m,n}$
if $1<m\le l$.  Finally, everything deforms to $d_1(0:1)$, which
has no deformations.

\subsubsection{Deformations of Codifferentials of the second kind}
\begin{figure}
\begin{center}
\epsfxsize=4.9in \epsfysize=3.0in \epsfbox{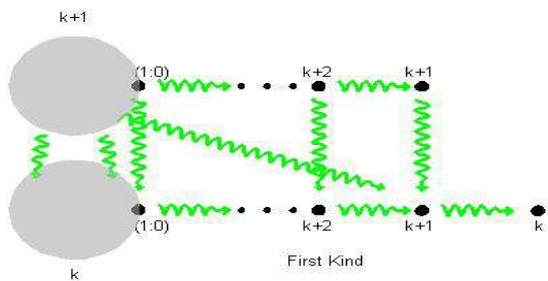}
\caption{{\bf The moduli space of codifferentials on} $\mathbf{\W_{-2}\oplus W_{-1}\oplus W_0}$.
The bottom row, depicting level $k$, contains
a sphere $d_k(\lambda:\mu)$ with marked point
$d_k(1:0)$, an infinite line of
dots $\{d_{k,m}|k\leq m\}$, and
horizontal arrows which show that $d_{k,l}$ deforms to
$d_{k,m}$ for $k<m<l$.
Although only one other row (depicting level $k+1$) is shown,
similar rows exist for all $k\in\mathbb N$.
Each arrow indicates the direction of a jump deformation.
For example, the arrows between the large spheres
signify that $d_{k+1}(\lambda:\mu)$ can deform to $d_k(\alpha:\beta)$, for
all values of $\alpha$  and $\beta$.
and the remaining downward and diagonal arrows represent that $d_{k,l}$
deforms to $d_{m,n}$ if $1<m<l$.
Although $d_{k,l}$ deforms to $d_m(\alpha:\beta)$ for all $1<m<k$,
the arrows are omitted in order to simplify the diagram.
Similarly, although $d_k(\lambda:\mu)$ can deform
to $d_m(\alpha:\beta)$ for all $1<m<k$, an infinite number of these
arrows between the spheres are omitted as well.
} \label{casetwo_label}
\end{center}
\end{figure}

A codifferential of the second kind is equivalent to one with
leading term of the form $d_k=\psa{1,k-1,0}3$. Since
\begin{align*}
H^n(d_k)_1=&\langle\psa{1,n-1,0}3\rangle,&n<k\\
H^n(d_k)_1=&0, &n\ge k,
\end{align*}
it follows that $d_k$ deforms to $d_n$ precisely when $n<k$.  Thus
the moduli space of codifferentials of the second kind is another
discrete chain, with jump deformations to lower elements on the
chain. A picture of the moduli space appears in figure
\ref{casetwo_label}.

\section{Conclusion}

Applications of infinity algebras occur in both \Z-grading and
\zt-grading, so it is important to understand the relationship
between them. We showed in this paper that even for 3-dimensional
infinity algebras this relationship is not trivial.
Sometimes  (as seen in section \ref{sec1bar2}) the natural map
between the moduli space of \Z-graded infinity algebras and the
moduli space of \zt-graded infinity algebras is an isomorphism.
In the example in section \ref{sec1bar2}, the moduli spaces are
quite complicated, but some seemingly unlikely coincidences
conspire to yield identical moduli spaces.  In fact, there are
twice as many \zt-graded codifferentials, but the group of
equivalences is just larger enough to identify them with the
\Z-graded codifferentials.

In general, the space of \Z-graded codifferentials will be much
sparser than the \zt-graded codifferentials, as is illustrated by
the example in section \ref{InjectiveExample}.  In that example,
the map from the moduli space of \Z-graded codifferentials to the
space of \zt-graded codifferentials was at least injective. When
the moduli space of \Z-graded codifferentials is sparse, it is not
really surprising that the map should be injective, but not
surjective.

In the example in section \ref{BadExample}, the map from the
moduli space of \Z-graded codifferentials fails to be injective in
a very interesting manner.  The map takes a stratum of the
\Z-graded moduli space which is a $\P^1$ to an orbifold stratum of
the \zt-graded moduli space given by a $\P^1/\Sigma_2$,  due to an
extra symmetry induced by the larger group of equivalences in the
\zt-graded case.  Thus even though the map fails to be injective,
the \zt-graded moduli space is closely related to the \Z-graded
moduli space.

Note that in each $3$-dimensional \linf\ algebra examined here, we
have shown that the moduli space is decomposable in a unique
manner as a union of orbifolds (cf. \cite{fp7}), with
the distinct strata glued together by jump deformations. In each
example, the map between the moduli spaces respects this
stratification, and in the one example where it fails to be
injective, the map has a very interesting structure on the stratum
on which it fails to be injective.

From the examples, we can deduce that the study of moduli spaces
of \Z-graded and \zt-graded moduli spaces cannot be reduced to a
study of only one of these types.  Nevertheless, it may be
possible to deduce important information about the \Z-graded
moduli space from that of the \zt-graded one.  Explorations of
higher dimensional examples may reveal more information about this
connection.

\section{Acknowledgements}

We would like to thank Jim Stasheff for reading various versions of this
manuscript and offering suggestions that greatly improved the readibility of the
paper. We also thank Carolyn Otto for creating the figures and helping with the
Maple programs.

\bibliographystyle{amsplain}

\begin{thebibliography}{10}

\bibitem{bfp1}
A.~Bodin, D.~Fialowski and M.~Penkava, \emph{Classification and versal
  deformations of \linf\ algebras on a $2|1$-dimensional space}, Homology,
  Homotopy and its Applications \textbf{7} (2005), no.~2, 55--86,
  math.QA/0401025.

\bibitem{dal2}
M.~Daily, \emph{\linf\ structures on spaces of low dimension}, Ph.D. thesis,
  North Carolina State University, 2004.

\bibitem{dal1}
\bysame, \emph{\linf\ structures on spaces with 3 one-dimensional components},
  Communications in Algebra \textbf{32} (2004), no.~5, 2041--2059,
  math.QA/0212030 v1.

\bibitem{fp1}
A.~Fialowski and M.~Penkava, \emph{Deformation theory of infinity algebras},
  Journal of Algebra \textbf{255} (2002), no.~1, 59--88, math.RT/0101097.


\bibitem{fp6}
\bysame, \emph{Strongly homotopy lie algebras of one even and two odd
  dimensions}, Journal of Algebra \textbf{283} (2005), 125--148,
  math.QA/0308016.

\bibitem{fp3}
\bysame, \emph{Versal deformations of three dimensional {L}ie algebras as
  \linf\ algebras}, Communications in Contemporary Mathematics \textbf{7}
  (2005), no.~2, 145--165, math.RT/0303346.


\bibitem{fp7}
\bysame, \emph{Examples of miniversal deformations of infinity algebras}, Forum
  Mathematica (2006), to appear, math RT/0510325.

\bibitem{lm}
T.~Lada and M.~Markl, \emph{Strongly homotopy {L}ie algebras}, Comm. in Algebra
  \textbf{23} (1995), 2147--2161.

\bibitem{ls}
T.~Lada and J.~Stasheff, \emph{Introduction to sh {L}ie algebras for
  physicists}, Int. J. Theor. Phys. \textbf{32} (1993), 1087--1103, Preprint
  hep-th 9209099.


\bibitem{ps2}
M.~Penkava and A.~Schwarz, \emph{\hbox{$A_\infty$} algebras and the cohomology
  of moduli spaces}, Dynkin Seminar, vol. 169, American Mathematical Society,
  1995, pp.~91--107.

\bibitem{sta4}
J.D. Stasheff, \emph{The intrinsic bracket on the deformation complex of an
  associative algebra}, Journal of Pure and Applied Algebra \textbf{89} (1993),
  231--235.

\end{thebibliography}
\providecommand{\bysame}{\leavevmode\hbox to3em{\hrulefill}\thinspace}
\providecommand{\MR}{\relax\ifhmode\unskip\space\fi MR }
\providecommand{\MRhref}[2]{%
  \href{http://www.ams.org/mathscinet-getitem?mr=#1}{#2}
}
\providecommand{\href}[2]{#2}

\end{document}